\documentclass{amsart}
\usepackage{amsmath}
\usepackage{latexsym}
\usepackage{amssymb}
\usepackage{amsthm}

\def\bg{\bmath\gamma}
\def\bmath#1{\mbox{\boldmath$#1$}}

\def\det{\operatorname{det}}

\def\e{\bmath \epsilon}
\def\ee{\mathbf e}

\def\t{{}^t\hskip-1pt}

\def\T{{}^t\hskip-2pt}

\def\Ta{\mathbf T}

\def\ve{\mathbf v}

\def\x{\mathbf x}
\def\y{\mathbf y}

\def\metric{I}  

\def\0{\mathbf 0}

\def\C{\mathbf C} 

\def\E{\mathbf E}

\def\GL{\mathbf{GL}}

\def\H{\mathbf H}
\def\Hom{\text{Hom}}

\def\Q{\mathcal{Q}}
\def\R{\mathbf R} 

\def\SO{\mathbf{SO}}
\def\s2{\mathbf S^2}

\def\sph{\mathbf S^3}
\def\Spin{\mathbf{Spin}}

\theoremstyle{plain}

\newtheorem*{theorem*}{Theorem}

\newtheorem{theorem}{Theorem}
\newtheorem{lemma}[theorem]{Lemma}
\newtheorem{corollary}[theorem]{Corollary}
\newtheorem{proposition}[theorem]{Proposition}
\newtheorem{definition}[theorem]{Definition}
\newtheorem{remark}[theorem]{Remark}

\begin{document}

\title{Compact surfaces with no Bonnet mate}


\author{Gary R. Jensen}
\address{(G. R. Jensen)
Department of Mathematics, Washington University, One Brookings Drive,
St. Louis, MO 63130, USA}
\email{gary@wustl.edu}

\author{Emilio Musso}
\address{(E. Musso) Dipartimento di Scienze Matematiche, Politecnico di Torino,
Corso Duca degli Abruzzi 24, I-10129 Torino, Italy}
\email{emilio.musso@polito.it}

\author{Lorenzo Nicolodi}
\address{(L. Nicolodi) Di\-par\-ti\-men\-to di Scienze Ma\-te\-ma\-ti\-che, Fisiche e Informatiche,
Uni\-ver\-si\-t\`a di Parma, Parco Area delle Scienze 53/A, I-43124 Parma, Italy}
\email{lorenzo.nicolodi@unipr.it}

\thanks{Authors partially supported by MIUR (Italy) under the PRIN project
\textit{Variet\`a reali e complesse: geometria, topologia e analisi armonica};
and by the GNSAGA of INdAM}



\begin{abstract} This note gives
  sufficient conditions (isothermic or totally nonisothermic) for an
  immersion of a
  compact surface to have no Bonnet mate.
\end{abstract}

\maketitle

\section{Introduction}

Consider a smooth immersion $\x:M \to \R^3$ of a connected, orientable
surface $M$, with unit normal vector field $\ee_3$.  Its induced metric
$\metric = d\x\cdot d\x$ and the orientation of $M$ induced by $\ee_3$ from
the standard orientation of $\R^3$ induce a complex structure on $M$,
which provides a decomposition into bidegrees of the second
fundamental form $II$ of $\x$ relative to $\ee_3$,
\[
-d\ee_3 \cdot d\x = II = II^{2,0} + H \metric + II^{0,2}.
\]
Here $H$ is the mean curvature of $\x$ relative to $\ee_3$ and
$II^{2,0} = \overline{II^{0,2}}$ is the \textit{Hopf quadratic differential} of
$\x$.
Relative to a complex chart $(U,z)$ in $M$,
\begin{equation}\label{eq:0}
\metric = e^{2u}dz d\bar z, \quad II^{2,0} = \frac12 he^{2u} dz dz,
\end{equation}
where the \textit{conformal factor} $e^u$, the \textit{Hopf invariant} $h$, and the mean
curvature $H$
satisfy the \textit{structure equations on $U$ relative to $z$},
\[
\aligned
-4e^{-2u}u_{z\bar z} &= H^2 - |h|^2 \;\;\text{ Gauss equation} \\
(e^{2u}h)_{\bar z} &= e^{2u}H_z \;\; \text{ Codazzi equation}
\endaligned
\]
and the Gauss curvature is $K = H^2 - |h|^2$.  See \cite[page 212]{MovingFrames}.

In 1867 Bonnet \cite{Bonnet1} began an investigation into the problem
of whether there exist noncongruent immersions $\x, \tilde \x:M \to
\R^3$ with the same induced metric, $\metric = \tilde\metric$, and the
same mean curvature, $H = \tilde H$.  This
\textit{Bonnet Problem} has been studied by Bianchi
\cite{Bianchi}, Graustein \cite{Graustein}, Cartan \cite{Cartan1942},
Lawson--Tribuzy \cite{LawsonTribuzy}, Chern \cite{Ch1},
Kamberov--Pedit--Pinkall \cite{KPP}, Bobenko--Eitner
\cite{BobenkoEitner,BobenkoEitner2}, Roussos--Hernandez \cite{RoussosHernandez},
Sabitov \cite{Sabitov2}, the present authors
\cite{MovingFrames}, and many others cited in these references.

\begin{definition}
An immersion $\x:M \to \R^3$ is \textit{Bonnet} if there is
a non\-con\-gru\-ent
immersion
$\tilde \x:M \to \R^3$ such that $\tilde \metric = \metric$
and $\tilde H = H$.  Then
$\tilde \x$ is called a
\textit{Bonnet mate} of $\x$ and $(\x,\tilde \x)$ form a \textit{Bonnet pair}.
\end{definition}

A constant mean curvature
(CMC) immersion $\x:M \to \R^3$, for which $M$ is simply connected and
$\x$ is not totally umbilic, admits
a 1-parameter
family of Bonnet mates, which are known as the associates of $\x$
\cite[Example 10.11, page 302]{MovingFrames}.
The local problem is thus to determine if an immersion $\x$ with
  nonconstant mean curvature has a Bonnet mate.
  By nonconstant mean curvature $H$ we mean that
  $dH \neq 0$ on a dense, open subset of $M$.

\begin{definition} A Bonnet immersion $\x:M \to \R^3$ is
  \textit{proper} if its mean curvature is nonconstant and there
  exist at least two noncongruent Bonnet mates.
\end{definition}

It is known \cite[page 211]{MovingFrames} that the umbilics of $\x$ are
precisely the zeros of its Hopf quadratic differential $II^{2,0}$.
For the following definitions we assume that $\x$ has no umbilics in
the domain $U$.  If
$(U,z)$ is a complex coordinate chart in $M$, then the local
coefficient $e^{2u}h$ of $2II^{2,0}$ in $U$ has the polar representation
\[
e^{2u}h = e^{G+ig},
\]
for a smooth function $G:U \to \R$ and a smooth map $e^{ig}:U \to
\mathbf S^1$.  The function $g:U \to \R$ is defined only locally, up
to an additive integral multiple of $2\pi$.  If $w= w(z)$ is another complex
coordinate in $U$, and if the invariants relative to it are denoted by
$\hat u$ and $\hat h$, then
\[
e^{2u}h = e^{2\hat u}\hat h (w')^2,
\]
where $w' = \frac{dw}{dz}$ is a nowhere zero holomorphic function of $z$.
Setting $e^{2\hat u}\hat h = e^{\hat G + i\hat g}$ on $U$,
we find by an elementary calculation
\begin{equation}\label{eq:0.1}
g_{\bar z z} = \hat g_{\bar z z}
\end{equation}
on $U$.  The Laplace-Beltrami operator of $(M,\metric)$ is given in
the local chart $(U,z)$ by $\Delta =
4e^{-2u}\frac{\partial^2}{\partial z\partial\bar z}$.  We conclude
from \eqref{eq:0.1} that $\Delta g = \Delta \hat g$ on $U$, and
therefore that $\Delta g$ is a globally defined smooth function on
$M$ away from the umbilic points of $\x$.

\begin{definition}\label{def0}
A surface immersion $\x:M \to \R^3$ is called \textit{isothermic} if
it has an atlas of charts $(U,(x,y))$ each of which satisfies $\metric
= e^{2u}(dx^2 + dy^2)$ and $II = e^u(adx^2 + cdy^2)$ \cite[Definition
  9.5, page 277]{MovingFrames}.
\end{definition}

Definition \ref{def0} is equivalent to the
following definition if there are no umbilics \cite[Corollary 9.14,
  page 280]{MovingFrames}.

\begin{definition}\label{def1}
  An umbilic free immersion $\x:M \to \R^3$ of an oriented connected
  surface is \textit{isothermic} if $\Delta g = 0$ identically on
  $M$.  It is \textit{totally nonisothermic} if $\Delta g \neq 0$ on
  a connected, open, dense subset of $M$.
\end{definition}

The following is known about umbilic free immersions $\x:M \to \R^3$
for which $M$
is simply connected. Cartan \cite{Cartan1942} proved that if $\x$ is
proper Bonnet, then it
has a 1-parameter family of
distinct mates \cite[Theorem 10.42, pages 340-342]{MovingFrames}.
Graustein \cite{Graustein} proved that if $\x$ is isothermic and
Bonnet, then it is proper Bonnet.  The present authors
\cite[Theorem 10.13, pages 303-304]{MovingFrames} proved that
if $\x$ is totally nonisothermic, then it has a unique Bonnet mate.

What is the global
situation?  In particular, if $M$ is compact, can an immersion $\x:M
\to \R^3$ have a Bonnet mate?  It is known, and proved in the next
section, that a necessary condition that $\x$ be Bonnet is that its
set of umbilics is a discrete subset of $M$.
Lawson--Tribuzy \cite{LawsonTribuzy}
proved that $\x$ cannot be proper Bonnet if $M$ is compact.
Roussos--Hernandez \cite{RoussosHernandez} proved that $\x:M \to \R^3$
has no
Bonnet mate if $M$ is compact and $\x$ is a surface of revolution with
nonconstant mean curvature.
Sabitov \cite[Theorem 13, page 144]{Sabitov2} gives a sufficient
condition preventing the existence of a  Bonnet mate when the mean
curvature is nonconstant and $M$ is
compact.  He gives no geometric interpretation of his condition.

The goal of this paper is to prove the following result.  It
generalizes the Roussos--Hernandez result, since a surface of
revolution is isothermic \cite[Example 9.7, page 277]{MovingFrames}.
It also gives a geometrical clarification
of the Sabitov result.

\begin{theorem*}
  Let $\x:M \to \R^3$ be a smooth immersion with nonconstant mean
  curvature $H$
  of a compact, connected surface, and suppose that $\mathcal D$, the set of
  umbilics of $\x$, is a discrete subset of $M$.
  \begin{enumerate}
    \item If $\x:M\setminus \mathcal D$ is isothermic, then $\x$ has no
      Bonnet mate.
    \item   If $\x$ is totally nonisothermic,
      then it has no Bonnet mate.
      \end{enumerate}
\end{theorem*}

\section{The deformation quadratic differential}

From the Gauss equation above, the Hopf
invariants $h$ and $\tilde h$ relative to a complex coordinate $z$
of two immersions with the same induced metric and the same
mean curvatures must satisfy
\[
|\tilde h| = |h|,
\]
since $\tilde u = u$.  Hence,
the only possible difference in the invariants of two such
immersions must be in
the arguments of the complex valued functions $h$ and $\tilde
h$.  Moreover, taking the difference of their Codazzi equations, we
get
\[
(e^{2u}\tilde h - e^{2u}h)_{\bar z} = e^{2u}(H_z - H_z) = 0,
\]
at every point of the domain $U$ of the complex coordinate $z$.  This
means that the function
\[
F = e^{2u}(\tilde h - h):U \to \C
\]
is holomorphic.

\begin{definition}\label{de:Bo:defquad}
If $\x, \tilde \x:M \to \R^3$
are immersions that induce the same complex structure on $M$,
then their \textit{deformation quadratic differential} is
\[
\mathcal Q = \widetilde{II}^{2,0} - II^{2,0}.
\]
\end{definition}

If $\x$ and $\tilde\x$ have the same induced metric and mean
curvature, then
the expression for $\mathcal Q$ relative to a complex
coordinate $z$ is
\begin{equation}\label{eq:Bo:0.1}
\mathcal Q = \frac12 e^{2u}(\tilde h - h)dzdz = \frac12 F dz dz,
\end{equation}
which shows that $\mathcal Q$
is a holomorphic quadratic differential on $M$, and
\begin{equation}\label{eq:Bo:0.2}
|F+e^{2u}h| = |e^{2u} \tilde h| = |e^{2u}h|
\end{equation}
on $U$, since $|\tilde h| = |h|$.
$\mathcal Q$ is identically zero on $M$ if and only if
$\tilde h = h$ in any complex coordinate system.
Therefore, by
Bonnet's Congruence Theorem, 
$\mathcal Q = 0$ if and only if
the immersions $\x$ and $\tilde \x$ are congruent in the sense that
there exists a rigid motion $(\y,A) \in \E(3)$ such that $\tilde \x =
\y + A\x:M \to \R^3$.  Thus, an immersion $\tilde \x:M \to \R^3$ is a
Bonnet mate of $\x:M \to \R^3$ if it induces the same
metric and mean curvature and the deformation quadratic differential
is not identically zero.

\begin{proposition}\label{pr:Bo:0.a} If an immersion $\x:M \to \R^3$
possesses a Bonnet mate $\tilde \x:M \to \R^3$, then the umbilics
of $\x$ must be isolated and coincide with those of $\tilde\x$.
\end{proposition}

\proof Under the given assumptions, the holomorphic quadratic
differential $\mathcal Q$ is not identically zero.  Therefore, in any
complex coordinate chart $(U,z)$, we have $\mathcal Q = \frac12 F dzdz$,
where
$F$ is a nonzero holomorphic function of $z$.  Its zeros
must be isolated.  A point $m \in U$ is an umbilic of $\x$ if and only if
$h(m)=0$ if and only if $\tilde h(m) = 0$, by~\eqref{eq:Bo:0.2}.  In either
case $F(m) = 0$ by~\eqref{eq:Bo:0.2}.  Therefore,
the set of umbilic points is a subset of the set of zeros of $\mathcal
Q$, which is a discrete subset of $M$.
\endproof

Let $\x:M \to \R^3$ be an immersion with a Bonnet mate $\tilde\x:M \to
\R^3$.  Let $(U,z)$ be a complex coordinate chart in $M$ and let
$h$ and $\tilde h$ be the Hopf invariants of $\x$ and $\tilde\x$, respectively,
relative to $z$ on $U$.  Let $\mathcal D$ be the set of umbilics of
$\x$, necessarily a discrete subset of $M$.
On $U \setminus \mathcal D$ we have $h$ never zero and
\[
\tilde h = hA,
\]
for a smooth function $A:U\setminus \mathcal D \to \mathbf S^1$,
where $\mathbf S^1 \subset \C$ is the unit circle.
On $U \setminus \mathcal D$ then, the difference of the Hopf
differentials is the holomorphic quadratic differential
\[
\mathcal Q =\widetilde{II^{2,0}} - II^{2,0} = II^{2,0}(A-1).
\]
This shows that $A:M \setminus \mathcal D \to \mathbf S^1$ is a
well-defined smooth map on all of $M \setminus \mathcal D$.

\begin{remark}\label{rem1}  Under our assumption of nonconstant $H$,
the map $A$ cannot
be constant, for otherwise $II^{2,0}$ would then be holomorphic and
thus $H$ would be constant by the Codazzi equation.
\end{remark}

\begin{proposition}[Sabitov\cite{Sabitov2}]\label{pr:Bo:2}
  If an immersion $\x:M \to \R^3$
possesses a Bonnet mate $\tilde \x:M \to \R^3$, then
the deformation quadratic differential $\mathcal Q$ of $\x$ is zero only
at the umbilics of $\x$.  Therefore, $A:M \setminus \mathcal D
\to \mathbf S^1$ never takes the value $1 \in \mathbf S^1$.
\end{proposition}

\proof
This is Theorem 1, pages 113ff of \cite{Sabitov2}.  He says the result
is stated in \cite{Bobenko3}, but he believes the proof there is inadequate.
Sabitov's proof uses results from the Hilbert boundary-value
problem. The following proof is essentially the same as Sabitov's, but
avoids use of the Hilbert boundary-value problem.

Seeking a contradiction,
suppose $\mathcal Q(m_0)=0$ for some point $m_0 \in M \setminus
\mathcal D$.  Since $\mathcal Q$ is holomorphic, and not identically zero, its
zeros are isolated.  Let $(U,z)$ be a complex coordinate chart of $M
\setminus \mathcal D$
centered at $m_0$, containing no other zeros of $\mathcal Q$, and such
that $z(U)$ is an open disk of $\C$.  Now $A(m_0) = 1$ and $A$ is
continuous, so we may assume $U$ chosen small enough that $A$ never
takes the value $-1$ on $U$.
Then there exists a smooth map
$v:U \to \R$ such that $-\pi <v<\pi$ and $A = e^{iv}$ on
$U$.  Since $A = 1$ on $U$ only at $m_0$, it follows that
\begin{equation}\label{eq3.1}
v(U\setminus\{m_0\}) \subset (-\pi,0) \text{ or } v(U\setminus\{m_0\})
\subset (0,\pi).
\end{equation}
Let $e^{2u}$ and $h$ be the conformal factor and
  Hopf invariant of $\x$ relative to $z$.  Then $h$ never zero on $U$
  implies it has a polar representation $h = e^{f+ig}$, for some
  smooth functions $f,g:U \to \R$.  Now $\mathcal Q =\frac12 Fdzdz$,
  where
  \[
  F = e^{2u}e^{f+ig}(e^{iv}-1)= e^{2u+f}(e^{i(g+v)} - e^{ig}):U \to \C
    \]
    is holomorphic.  Using the identity
    \[
    e^{i(g+v)} - e^{ig} = e^{i(2g+v)/2}(e^{iv/2} - e^{-iv/2}) =
      2ie^{i(g+v/2)}\sin(v/2),
      \]
      we get
      \[
      F = 2ie^{2u+f+i(g+v/2)}\sin(v/2)
      \]
      on $U$.  The contour integral of
      $d\log F$ about any circle in $U$ centered at $m_0$ is $2\pi i$ times
      the number of zeros of $F$ inside the circle.  By assumption,
      this integral is not zero.  But,
      \[
      d\log F = d(2u+f+i(g+v/2)) + d\log(|\sin(v/2)|),
      \]
      and the contour integral of the right hand side is zero, since
      these are exact differentials on $U \setminus\{m_0\}$.  In fact,
      the values of $v/2$ on $U \setminus \{m_0\}$ lie entirely in
      $(0,\pi/2)$ or entirely in $(-\pi/2,0)$, so $\sin(v/2)$ is never zero.
      This is
      the desired contradiction to our assumption that $\mathcal Q$
      has a zero in $M \setminus \mathcal D$.
\endproof

As a consequence of this Proposition, the smooth map $A:M
\setminus \mathcal D \to \mathbf S^1$ never takes the value $1 \in
\mathbf S^1$, so there exists a smooth map
\[
r:M\setminus \mathcal D \to (0,2\pi) \subset \R,
\]
such that $A = e^{ir}$ on $M\setminus \mathcal D$.

\section{Proof of the Theorem}

\proof Seeking a contradiction, we
suppose that $\x$ possesses a Bonnet mate $\tilde\x:M \to \R^3$.
Let $II^{2,0}$ and $\widetilde{II^{2,0}}$ be the Hopf quadratic
differentials of $\x$ and $\tilde \x$, respectively.  By the preceding
propositions, the quadratic differential $\widetilde{II^{2,0}} -
II^{2,0}$ is holomorphic on $M$, and on $M \setminus \mathcal D$
\[
\widetilde{II^{2,0}} - II^{2,0} = II^{2,0}(e^{ir}-1),
\]
where the
function $r:M\setminus \mathcal D \to (0,2\pi)$ is smooth.  Let
$(U,z)$ be a complex coordinate chart in $M\setminus \mathcal D$.
Let $h$ and $e^u$ be the Hopf invariant and conformal factor of $\x$ relative
to $z$.
Then $h = e^{f+ig}$ on $U$, for some
smooth functions $f:U \to \R$ and
$e^{ig}:U \to \mathbf S^1$.

1). If $\x$ is isothermic, then $g_{\bar z z}=0$ identically on $U$.
Let $G = f+2u:U \to \R$.
Then $(e^{G+ig}(e^{ir}-1))_{\bar z} = 0$ implies
\begin{equation}\label{eq1}
  r_{\bar z} = i(G+ig)_{\bar z}(1-e^{-ir})
\end{equation}
on $U$.
Applying $\partial_z$ to this, and using that $r_z$ is the complex
conjugate of $r_{\bar z}$, we find
\begin{equation}\label{eq1.a}
  r_{\bar z z} = 0
\end{equation}
on $U$.  Hence, $r:M \setminus \mathcal D \to (0,2\pi)$ is a bounded
harmonic function.  Since the points of $\mathcal D$ are isolated and
$r$ is bounded, we know that $r$ extends to a harmonic function on all
of $M$.  But then $r$ must be constant, since $M$ is compact.  This
contradicts our assumption of nonconstant $H$, by Remark \ref{rem1}.

2).
If $\x$ is totally nonisothermic, we have either $\Delta g \leq 0$ or $\Delta
g \geq 0$
on $M\setminus \mathcal D$.  To be specific, let us suppose
that $\Delta g \leq 0$ on $M\setminus \mathcal D$.
Now
\eqref{eq1} holds and by the proof of Theorem
10.13 on pages 303-304 of \cite{MovingFrames}, we have
\begin{equation}\label{eq2}
  e^{ir} = 1 + \frac{-2 g_{\bar z z}}{D}(g_{\bar z z} + iL),
  \end{equation}
on $U$, where
$L = |G_{\bar z}+ ig_{\bar z} |^2 - G_{\bar z z}$ and $D = g_{\bar z
  z}^2 + L^2$.
Applying
$\partial_z$ to \eqref{eq1} and using \eqref{eq2}, we find
\begin{equation}\label{eq3}
  r_{\bar z z} = -2g_{\bar z z},
\end{equation}
on $U$.  Therefore, $\Delta r = -2\Delta g \geq 0$ on $M\setminus
\mathcal D$.

Recall \cite[Def. \S2.1, pages 40-41]{HaymanKennedy} that a function
$v:V\to \R\cup \{-\infty\}$ on a domain $V \subset \C$ is
\textit{subharmonic} if
\begin{enumerate}
\item $-\infty \leq v(z) < +\infty$ in $V$.
\item $v$ is upper semi-continuous in $V$. (This means that
  for any $c \in \R$, the set $\{z \in U: v(z)<c\}$ is open in V.)

  \item If $z_0$ is any point of $V$ then there exist arbitrarily
    small positive values of $R$ such that
    \[
    v(z_0) \leq \frac1{2\pi R} \int_0^{2\pi} v(z_0 + Re^{it}) dt.
    \]
\end{enumerate}
If $v$ is of class $C^2$ in $V$, then $v$ is subharmonic in $V$ if and only
if $v_{\bar z z} \geq 0$ in $V$ \cite[Example 3, page
  41]{HaymanKennedy}.

If $M$ is a connected Riemann surface, we define a function $v:M \to
\R\cup \{-\infty\}$ to be subharmonic if for any complex coordinate chart
$(U,z)$ of $M$, the local representative $v\circ z^{-1}:z(U)\to\R$ is
subharmonic. This is well-defined by the Corollary to Theorem 2.8 on page 53
of \cite{HaymanKennedy}.

We conclude from \eqref{eq3} that $r$ is subharmonic on
$M\setminus \mathcal D$.
In the event that $\Delta g \geq 0$ on $M\setminus \mathcal D$, we
conclude that $-r$ is subharmonic and continue as below with $-r$.

Suppose $(U,z)$ is a complex
coordinate chart centered at a point $m_0 \in \mathcal D$, and small
enough that no other point of $\mathcal D$ lies in it.  Then $r\circ
z^{-1}$ is subharmonic on the open set $z(U) \setminus\{0\}$, so it extends
uniquely to a subharmonic function on $z(U)$, by Theorem 5.8 on page 237 of
\cite{HaymanKennedy}.  It follows that $r$ extends uniquely to a
subharmonic function on $M$.

By Theorem 1.2 on page 4 of \cite{HaymanKennedy}, if $v:V \to \R \cup
\{-\infty \}$ is
upper semi-continuous on a nonempty compact domain $V \subset \C$,
then $v$ attains
its maximum on $V$; i.e., there exists $z_0\in V$ such that $v(z) \leq
v(z_0)$ for all $z \in V$.  The same proof shows that this is true for
an upper semi-continuous function on a compact Riemann surface.  Thus,
the subharmonic
function $r:M \to \R\cup \{-\infty\}$ attains its maximum at some
point $m_0 \in M$.
Let $(U,z)$ be a complex coordinate chart centered at $m_0$.  Choose
$R>0$ such that the disk $D(0,R) = \{z\in \C: |z|\leq R\}$ is
contained in $z(U)$.  By the maximum principle for subharmonic functions
\cite[Theorem 2.3, page 47]{HaymanKennedy}, $r\circ z^{-1}$ must be
constantly equal to $r(m_0)$ on $D(0,R)$.  It follows that
\[
E =\{m\in M: r(m) = r(m_0)\}
\]
is an open subset of $M$.  But
\[
E = M \setminus \{m \in M: r(m) < r(m_0)\}
\]
is closed, since $r$ is upper semi-continuous.  We conclude that $r$
is constant on
$M$, which is our sought for contradiction, by Remark \ref{rem1}.

\endproof

\bibliography{Bibliography}
\bibliographystyle{alpha}

\end{document}